\documentclass[final,1p,times]{elsarticle}

\usepackage{amssymb}
\usepackage{amsmath}
\usepackage{xcolor}

\journal{Numerical Functional Analysis and Optimization}

\newcommand{\N}{\mathbb{N}}
\newcommand{\R}{\mathbb{R}}
\newcommand{\be}{\begin{equation}}
\newcommand{\ee}{\end{equation}}

\begin{document}

\newtheorem{thm}{Theorem}[section]
\newtheorem{lem}{Lemma}[section]
\newtheorem{rmk}{Remark}[section]
\newtheorem{pro}{Proposition}[section]
\newtheorem{example}{Example}[section]
\newdefinition{definition}{Definition}[section]
\newproof{pf}{Proof}
\renewcommand{\theequation}{\thesection.\arabic{equation}}

\begin{frontmatter}



\title{Local average sampling and reconstruction with fundamental splines of fractional order}


\author[dev]{P. Devaraj}
\ead{devarajp@iisertvm.ac.in}
\author[mas]{P. Massopust}
\ead{massopust@ma.tum.de}
\author[yug]{S. Yugesh\corref{}}
\ead{mathsyugesh@gmail.com, yugeshs@ssn.edu.in}

\cortext[]{Corresponding Author}
\address[dev]{School of Mathematics, Indian Institute of Science Education and Research Thiruvananthapuram,Vithura, Thiruvananthapuram-695551.}
\address[mas]{
Center of Mathematics, Research Unit M15, Technical University Munich, Boltzmannstrasse 3,
85748 Garching b. Munich, Germany}
\address[yug]{Department of Mathematics, SSN College of Engineering, Kalavakkam-603 110, Tamil Nadu, India.}

\begin{abstract}
We analyse sampling and average sampling techniques for fractional spline subspaces of $L^{2}({\mathbb{R}}).$ Fractional B-splines $\beta_{\sigma}$ are  extensions of Schoenberg's polynomial splines of integral order to real order $\sigma > -1$. We present the interpolation with fundamental splines of fractional order for $\sigma \geq 1$ and the average sampling with fundamental splines of fractional order for $\sigma \geq \frac{3}{2}.$ Further, we generalise Kramer's lemma in the context of  local average sampling.
\end{abstract}


\end{frontmatter}

\section{Introduction and Preliminaries}
In digital signal and image processing, continuous signals need to be represented by their discrete samples. A fundamental problem  is how to represent a continuous signal in terms of its discrete samples. The main goal of sampling theory is to reconstruct a function from a suitable class of functions  from its discrete samples \cite{gro, buz}.  The well-known Shannon sampling theorem states that any band-limited signal $f$ is completely determined by its samples \cite{ald2, ald3, uns1}. In practical situations, the available signals need not be band-limited. In order to handle such situations, many authors have discussed the sampling and reconstruction problem in general shift invariant spaces and spline subspaces \cite{ald, ald1, ald2, ald3, Ald4, Ald5, uns, uns1, gro, buz}. The shift invariant spaces and spline spaces yield many advantages in practical applications. 
 
In \cite{Ald4}, Aldroubi et al. have studied the problem of reconstructing  functions from a set of nonuniformly distributed weighted average samples in the context of shift invariant subspaces of $L^{p}({\mathbb{R}}^{d})$ generated by $p$-frames.  Moreover, they have developed fast approximation-projection iterative reconstruction algorithms. 
    The same authors have analysed, in \cite{Ald5}, uniformly sampled convolution and stable average samplers and their reconstructions  over shift invariant spaces. Furthermore, they also studied sampling and reconstruction on irregular grids  and  established the connections between stable deconvolution and stable reconstruction from samples  after convolution is subtle.
    
    In \cite{sch}, Schoenberg introduced cardinal polynomial splines which are compactly supported functions. Spline functions are a very convenient tool for solving practical application problems. Many  studies were done on sampling and reconstruction theorem for spline subspaces \cite{ald1, ald3, uns, uns1, gro, buz, sun, sun1, chu}.

\vskip 1em

In general, the degrees of splines are integers. An extension of polynomial B-splines to fractional order is known as fractional B-splines. M. Unser and Th. Blu introduced such fractional B-splines in \cite{uns2}. They showed that all of desirable computational properties of cardinal B-splines of integral degree carry over to the fractional B-splines. But in general, fractional B-splines are not compactly supported functions.

Symmetric fractional B-splines are defined in the Fourier domain in \cite{uns2} by
\[
\widehat{\beta}_{\sigma}(\omega):= \left|\frac{\sin(\frac{\omega}{2})}{\frac{\omega}{2}}\right|^{\sigma+1},
\]
for $\sigma > -1$ and it is shown that  $\beta_{\sigma} \in L^{1}({\mathbb{R}}) \cap L^{2}({\mathbb{R}}).$ In the following, we consider the fractional spline space $V$ which is a subspace of $L^{2}({\mathbb{R}})$ with a generator $\beta_{\sigma} \in L^{2}({\mathbb{R}}),$
$$
V:= \left\{f(t)= \sum\limits_{n \in {\mathbb{Z}}} a_{n}\beta_{\sigma}(t-n) : \{a_{n}\} \in l^{2}({\mathbb{Z}}) \right\},\qquad \sigma > -1.
$$
For this range of $\sigma$-values, this representation is stable and the fractional spline space $V$ is a  well-defined subspace of $L^{2}({\mathbb{R}})$ \cite{uns2}.

%

\noindent
\section{Interpolation of fractional order fundamental splines}
We consider the interpolation with fundamental splines of fractional order. In \cite{for}, the interpolation with fundamental splines of fractional order of the form
\be\label{B}
\widehat{B}_\sigma(\omega ) := \int_{\R} B_\sigma(t)e^{-i\omega t}\, dt := \left( \frac{1-e^{-i\omega}}{i\omega}\right)^\sigma = e^{-i \sigma \pi/2}\,\left(\frac{\sin(\omega/2)}{\omega/2}\right)^\sigma,
\ee
is analysed for $\sigma \geq 2$ and $\sigma \notin 2{\mathbb{N}}+1.$ The fractional splines considered in the present work are symmetric about the $y$-axis and those considered  in \cite{for} are defined on $[0,\infty).$ In this paper, we analyse the analogues results  for $\sigma \geq 1.$

\vskip 1em

 The main aim in interpolation is to construct a fundamental cardinal spline of fractional order $L_{\sigma}:{\mathbb{R}}\rightarrow {\mathbb{R}},$
\begin{equation}
L_{\sigma}:= \sum\limits_{k \in {\mathbb{Z}}} c_{k}^{(\sigma)}\beta_{\sigma}(.-k),
\label{eqn1}
\end{equation}
satisfying the interpolation condition:
\begin{equation}
L_{\sigma}(n)=\delta_{n,0}, \mbox{   } n \in {\mathbb{Z}},
\label{finter}
\end{equation}
for an appropriate bi-infinite sequence $\{c_{k}^{(\sigma)}: \mbox{   } k \in {\mathbb{Z}}\}$ and for
suitable values of $\sigma.$

\indent Considering the formal series (\ref{eqn1}) and taking the Fourier transform on both sides of equation (\ref{eqn1}), we get
\begin{equation}
\widehat{L}_{\sigma}(\omega)=\widehat{\beta}_{\sigma}(\omega) \sum\limits_{k \in {\mathbb{Z}}} c_{k}^{(\sigma)}e^{-ik \omega}.
\label{eqn3}
\end{equation}
In order to find the fundamental splines, we obtain from   the formal series (\ref{eqn1}) and (\ref{finter}),
$$\sum\limits_{k \in {\mathbb{Z}}} c_{k}^{(\sigma)}e^{-ik \omega}=\frac{1}{\sum\limits_{k \in {\mathbb{Z}}}\beta_{\sigma}(k)e^{-ik \omega}}.$$
i.e., on the unit circle $|z|=1,$ we must have
$$\sum\limits_{k \in {\mathbb{Z}}} c_{k}^{(\sigma)}z^{-k}=\frac{1}{\sum\limits_{k \in {\mathbb{Z}}}\beta_{\sigma}(k)z^{-k}}.
$$
Therefore using equation (\ref{eqn3}), we obtain
\begin{equation}
\widehat{L}_{\sigma}(\omega)=\frac{\widehat{\beta}_{\sigma}(\omega)}{\sum\limits_{k=-\infty}^{\infty}\beta_{\sigma}(k)z^{-k}}.
\label{eqn4}
\end{equation}
In order to show that $L_{\sigma}$ is well defined we have to show that the denominator in (\ref{eqn4}) is not zero on the unit circle $|z|=1$. We obtain such a sufficient condition in the following theorem.

%

\begin{thm}
Let  $\sigma \geq 1$. Then $G(z) := \sum\limits_{n \in {\mathbb{Z}}} \beta_{\sigma}(n)z^{-n}$ has no roots on the unit circle $|z|=1.$
\label{thm2}
\end{thm}

\begin{pf}

 In terms of the Fourier transforms we can write $$\displaystyle{\widehat{L}_{\sigma}(\omega)=\frac{\widehat{\beta}_{\sigma}(\omega)}{\sum\limits_{k \in {\mathbb{Z}}} \widehat{\beta}_{\sigma}(\omega+2\pi k)}= \frac{\frac{1}{|\omega|^{\sigma+1}}}{\sum\limits_{k \in {\mathbb{Z}}}\frac{1}{|\omega+2\pi k|^{\sigma+1}}}.}$$
\noindent
Here, using the same arguments as in \cite{for}, it suffices to consider $0<\omega<2\pi.$

The denominator on the right can be written as
$$
\sum\limits_{k \in {\mathbb{Z}}}\frac{1}{|\omega+2\pi k|^{\sigma+1}}=\frac{1}{(2\pi)^{\sigma+1}}\sum\limits_{k \in {\mathbb{Z}}}\frac{1}{|a+ k|^{\sigma+1}},\;\;\;a :=\frac{\omega}{2\pi}. $$

\noindent Clearly $0<a<1.$ The latter sum has to be zero free for all $0<a<1.$

The following manipulations hold true:
\begin{eqnarray*}
\sum\limits_{k \in {\mathbb{Z}}}\frac{1}{|a+ k|^{\sigma+1}}&=&\sum\limits_{k =-\infty}^{-1}\frac{1}{|a+ k|^{\sigma+1}} + \sum\limits_{k=0}^{\infty}\frac{1}{(a+ k)^{\sigma+1}}\\
 &=& \sum\limits_{k=1}^{\infty}\frac{1}{(k-a)^{\sigma+1}}+\sum\limits_{k=0}^{\infty}\frac{1}{(a+ k)^{\sigma+1}}\\
 &=& \sum\limits_{k=0}^{\infty}\frac{1}{(k+1-a)^{\sigma+1}}+ \sum\limits_{k=0}^{\infty}\frac{1}{(a+ k)^{\sigma+1}} \\
 &=&\zeta(\sigma +1,1-a)+ \zeta(\sigma+1,a) =: Z(\sigma + 1, 1-a,a),
\end{eqnarray*}
where $\zeta(\sigma+1,a)$ denotes the classical Hurwitz zeta function \cite{baer}.
Using the result stated in Lemma 2 of \cite{for}, we see that if $\sigma+1\ge 1+a$ and  $\sigma +1 \ge 1+ 1-a$ for all $0<a<1,$ i.e., $\sigma\ge 1,$ then both zeta functions are zero free for  $\sigma \ge 1.$
\noindent
The arguments employed in \cite{for} to show that $Z(\sigma + 1, 1-a,a)$ is zero free for $0 < a < \frac12$ and for $\frac12 < a <1$ also apply to the current setting. For $a=\frac{1}{2}$, one obtains $Z(\sigma+1,\frac{1}{2},\frac{1}{2})=2\ \zeta(\sigma+1,\frac{1}{2})\neq 0 $ for $\sigma\ge 1.$
\hfill{$\square$}
\end{pf}

\begin{rmk}
Note that unlike in \cite{for}, there is no further restriction on $\sigma$. The factor $e^{-i \sigma \pi/2}$ in the definition of $\widehat{B}_\sigma$, is responsible for the exclusion of the odd integer powers $\sigma = 2n + 1$, $n\in \N$.
\end{rmk}

The fundamental cardinal spline $L_{\sigma}$ of fractional order $\sigma$ with $\sigma \geq 1$ is an element of $L^{1}({\mathbb{R}}) \cap L^{2}({\mathbb{R}}),$ since $\beta_{\sigma} \in L^{1}({\mathbb{R}}) \cap L^{2}({\mathbb{R}}).$ Further, $L_{\sigma}$ is uniformly continuous on ${\mathbb{R}}.$

\vskip 1em

\indent In order to  analyse the sampling theorem for the fundamental cardinal spline $L_{\sigma},$  we consider the following version of Kramer's lemma \cite{kra} which appears in \cite{for,gar1}.
\begin{thm}
Let $\emptyset \neq I,$ $\Omega \subseteq {\mathbb{R}}$ and let $\{\phi_{k}: \mbox{  } k \in {\mathbb{Z}}\}$ be an orthonormal basis of $L^{2}(I).$ Suppose that $\{S_{k}: \mbox{  } k \in {\mathbb{Z}}\}$ is a sequence of functions $S_{k}:\Omega\rightarrow {\mathbb{C}}$ and ${\bf t}:=\{t_{k} \in {\mathbb{R}}: \mbox{  } k \in {\mathbb{Z}}\}$ a numerical sequence in $\Omega$ satisfying the conditions
\begin{itemize}
\item[\emph{C1.}] $S_{k}(t_{l})=a_{k}\delta_{kl},$ $(k,l) \in {\mathbb{Z}} \times {\mathbb{Z}},$ where $a_{k} \neq 0$;
\item[\emph{C2.}] ${\sum\limits_{k \in {\mathbb{Z}}}|S_{k}(t)|^{2}< \infty},$ for each $t \in \Omega.$
\end{itemize}
Define a function $K:I\times \Omega \rightarrow \mathbb{C}$ by
$$\displaystyle{K(x,t):=\sum\limits_{k \in {\mathbb{Z}}}S_{k}(t)\overline{\phi_{k}(x)}},$$
and a linear integral transform $\mathcal{K}$ on $L^{2}(I)$ by
$$(\mathcal{K}f)(t):= \int_{I}f(x)K(x,t)dx.$$
Then $\mathcal{K}$ is well defined and injective. Furthermore, if the range of $\mathcal{K}$ is denoted by
$$\mathcal{H}:=\{g: \mbox{  } {\mathbb{R}}\rightarrow {\mathbb{C}}: \mbox{  } g= \mathcal{K}f, \mbox{  } f \in L^{2}(I)\},$$
then
\begin{enumerate}
\item[\emph{(i)}]
$(\mathcal{H}, \langle \cdot, \cdot \rangle _{\mathcal{H}})$ is a Hilbert space and is isometrically isomorphic to $L^{2}(I),$ i.e, $\mathcal{H}\cong L^{2}(I),$ when endowed with the inner product
$$\langle F,G\rangle_{\mathcal{H}}:= \langle f,g \rangle_{L^{2}(I)},$$
where $F:= \mathcal{K}f$ and $G=\mathcal{K}g.$
\item[\emph{(ii)}]
$\{S_{k}: \mbox{  } k \in {\mathbb{Z}}\}$ is an orthonormal basis for $\mathcal{H}.$
\item[\emph{(iii)}]
Each function $ f \in \mathcal{H}$ can be recovered from its samples on the sequence $\{t_{k}: k \in {\mathbb{Z}}\}$ via the formula
$$\displaystyle{f(t)= \sum\limits_{k \in {\mathbb{Z}}}f(t_{k}) \frac{S_{k}(t)}{a_{k}}}.$$
The above series converges absolutely and uniformly on subsets of ${\mathbb{R}},$ where $||K(\cdot,t)||_{L^{2}(I)}$ is bounded.
\end{enumerate}
\label{samplingk}
\end{thm}
\indent If we take $\Omega:={\mathbb{R}},$ ${\bf t}:= {\mathbb{Z}},$ $a_{k}=1$ for all $k \in {\mathbb{Z}},$ and the interpolating functions $S_{k}=L_{\sigma}(\cdot -k),$ $k \in {\mathbb{Z}},$ then Theorem \ref{samplingk} implies the following theorem.
\begin{thm}
Let $\emptyset \neq I \subseteq {\mathbb{R}}$ and let $\{\phi_{k}: \mbox{  } k \in {\mathbb{Z}}\}$ be an orthonormal basis of $L^{2}(I).$ Let $L_{\sigma}$ denote the fundamental cardinal spline of fractional order $\sigma \geq 1$. Then the following hold:
\begin{enumerate}
\item[\emph{(i)}]
The family $\{L_{\sigma}(\cdot -k): k \in {\mathbb{Z}}\}$ is an orthonormal basis of the Hilbert space $(\mathcal{H},\langle \cdot,\cdot\rangle_{\mathcal{H}}),$ where $\mathcal{H} = \mathcal{K}(L^{2}(I))$ and $\mathcal{K}$ is  the injective integral operator
$$\displaystyle{\mathcal{K}f=\sum\limits_{k \in {\mathbb{Z}}}\langle f,\phi_{k}\rangle_{L^{2}(I)}L_{\sigma}(\cdot -k)}, \mbox{     } f \in L^{2}(I).$$
\item[\emph{(ii)}]
Every function $f \in \mathcal{H}\cong L^{2}(I)$ can be recovered from its samples on the integers via
\begin{equation}
f(\cdot)=\sum\limits_{k \in {\mathbb{Z}}}f(k)L_{\sigma}(\cdot -k),
\label{krameq}
\end{equation}
where the series (\ref{krameq}) converges absolutely and uniformly on all subsets of ${\mathbb{R}}.$
\end{enumerate}
\end{thm}
\begin{pf}
   Taking $S_{k}=L_{\sigma}(\cdot -k),$ $k \in {\mathbb{Z}},$  condition C1 in Theorem \ref{samplingk} is reduced to
$$L_{\sigma}(n)=\delta_{n,0}, \mbox{   } n \in {\mathbb{Z}}.$$
This condition is verified by the equation (\ref{finter}). Since the fundamental cardinal spline $L_{\sigma}$ is an element of $L^{1}({\mathbb{R}}) \cap L^{2}({\mathbb{R}}),$ the condition C2 is also established. Further, as the unfiltered splines $\{\beta_{\sigma}(\cdot -k): k \in {\mathbb{Z}}\}$ form a Riesz basis of the space $V$ (see, \cite{uns2}), $||K(\cdot,t)||_{L^{2}(I)}$ is bounded on ${\mathbb{R}}.$ Hence (i) and (ii) follow from Theorem \ref{samplingk}.
\hfill{$\square$}
\end{pf}

\section{Local average sampling for fractional spline space}
In practical situations, it is difficult to measure the exact values of the samples. The  measurement process depends on the aperture device used  for capturing the samples.  An appropriate model is to assume that the samples are local average samples of the form
$$f \star h(n) = \int_{-\infty}^{\infty} f(t)h(n-t)dt, \mbox{  } n \in {\mathbb{Z}},$$
where the averaging function $h(t)$ reflects the characteristics of the acquisition device. In this section we carry over the interpolation with fundamental splines of fractional order to the local average sampling context. We assume that the averaging function $h(t)$ is compactly supported in $L^{1}({\mathbb{R}}).$

\vskip 1em

 The fundamental spline of fractional order $L_{h,\sigma}: {\mathbb{R}}\longrightarrow {\mathbb{R}}$ for the average sampling problem is
\begin{equation}
L_{h,\sigma}:= \sum\limits_{k \in {\mathbb{Z}}} c_{k}^{(h,\sigma)}\beta_{\sigma}(\cdot -k)
\label{eqn5.10}
\end{equation}
satisfying the weighted interpolation condition:
\begin{equation}
L_{h,\sigma}\star h(n)=\delta_{n,0}, \mbox{   } n \in {\mathbb{Z}},
\label{finterh}
\end{equation}
Taking the  Fourier transform on both sides of the formal series (\ref{eqn5.10}), we obtain
\begin{equation}
\widehat{L}_{h,\sigma}(\omega)=\widehat{\beta}_{\sigma}(\omega) \sum\limits_{k \in {\mathbb{Z}}} c_{k}^{(h,\sigma)}e^{-ik \omega}.
\label{eqn5.11}
\end{equation}
In view of formal series (\ref{eqn5.10}) and (\ref{finterh}), we obtain
$$\sum\limits_{k \in {\mathbb{Z}}} c_{k}^{(h,\sigma)}e^{-ik \omega}=\frac{1}{\sum\limits_{k \in {\mathbb{Z}}}\beta_{\sigma} \star h(k)e^{-ik \omega}}$$
Therefore on the unit circle $|z|=1,$
$$
\sum\limits_{k \in {\mathbb{Z}}} c_{k}^{(h,\sigma)}z^{-k}=\frac{1}{\sum\limits_{k \in {\mathbb{Z}}}\beta_{\sigma} \star h(k)z^{-k}}.$$

Using equation (\ref{eqn5.11}) we obtain
\begin{equation}
\widehat{L}_{h,\sigma}(\omega)=\frac{\widehat{\beta}_{\sigma}(\omega)}{\sum\limits_{k \in {\mathbb{Z}}}\beta_{\sigma} \star h(k)z^{-k}}.
\label{eqn5.12}
\end{equation}
In order to construct the fundamental splines of fractional order for the local average sampling problem we have to   show that the denominator of (\ref{eqn5.12}) is not zero on the unit circle $|z|=1.$ We obtain sufficient conditions on  $\sigma$ and $h$ for which this holds in the following theorem.



\begin{thm}
Let $\sigma > 1$ and $h \in L^{1}({\mathbb{R}}).$ Consider the non-negative averaging function $h(t)$, whose support is contained in $[-l,l]$, where
\begin{equation}
l=\frac{1}{\pi}\arccos\left[\left(\frac{1}{2^{\sigma+1}}\right)^{\frac{1}{3}}\right].
\label{eqn3.2}
\end{equation}
Then $G_{h}(z)= \sum\limits_{n \in {\mathbb{Z}}} \beta_{\sigma} \star h(n)\, z^{-n}$ has no roots on the unit circle $|z|=1.$
\label{thm4}
\end{thm}

\begin{pf}

We consider  $F:=\beta_{\sigma} \star h.$  As $\beta_{\sigma}, h \in L^{1}({\mathbb{R}}),$ we get
$F \in L^{1}({\mathbb{R}})$
and $\widehat{F}\in L^{1}({\mathbb{R}}).$ Now the $2\pi$-periodic function defined by
\begin{eqnarray*}
\Psi_{\widehat{F}}(x)&=&\sum\limits_{k=-\infty}^{k=\infty}\widehat{F}(x+2\pi k)\\
&=&\sum\limits_{k=-\infty}^{k=\infty}\widehat{\beta}_{\sigma}(x+2\pi
k)\widehat{h}(x+2\pi k)
\end{eqnarray*}
 converges everywhere and its corresponding
Fourier series is given by
\begin{eqnarray*}
\Psi_{\widehat{F}}(x)&\sim&\sum\limits_{k=-\infty}^{k=\infty}c_{k}(\Psi_{\widehat{F}})e^{ikx},
\end{eqnarray*}
where
\begin{eqnarray*}
c_{k}(\Psi_{\widehat{F}})&=&\frac{1}{2\pi}
\int_{0}^{2\pi}e^{-ikx}\Psi_{\widehat{F}}(x)dx\\
&=&\frac{1}{2\pi}
\int_{0}^{2\pi}e^{-ikx}\sum\limits_{j=-\infty}^{j=\infty}\widehat{F}(x+2\pi
j)dx\\
&=&\frac{1}{2\pi}\sum\limits_{j=-\infty}^{j=\infty}\int_{0}^{2\pi}e^{-ikx}\widehat{F}(x+2\pi
j)dx\\
&=&\frac{1}{2\pi}\sum\limits_{j=-\infty}^{j=\infty}\int_{2\pi j}^{2\pi
(j+1) }e^{-ikx}\widehat{F}(x)dx\\
&=&\frac{1}{2\pi} \int_{-\infty}^{\infty }e^{-ikx}\widehat{F}(x)dx\\
&=&\frac{1}{2\pi}\widehat{\widehat{F}}(k).\\
\end{eqnarray*}
Therefore we obtain,
\begin{eqnarray*}
\Psi_{\widehat{F}}(x)&=&\sum\limits_{k=-\infty}^{k=\infty}c_{k}(\Psi_{\widehat{F}})e^{ikx}\\
&=&\sum\limits_{k=-\infty}^{k=\infty}\frac{1}{2\pi}\widehat{\widehat{F}}(k)e^{ikx}\\
&=&\sum\limits_{k=-\infty}^{k=\infty}F(k)e^{-ikx}.
\end{eqnarray*}
As the  Fourier series  converges everywhere, we get
\begin{eqnarray*}
\sum\limits_{k=-\infty}^{k=\infty}\widehat{\beta}_{\sigma}(x+2\pi
k)\widehat{h}(x+2\pi
k)&=&\sum\limits_{k=-\infty}^{k=\infty}F(k)e^{-ikx}\\
&=&\sum\limits_{k=-\infty}^{k=\infty}(\beta_{\sigma} \star h)(k)e^{-ikx}
\end{eqnarray*}
hence  we get
\begin{equation}
G_{h}(z)=\sum\limits_{k=-\infty}^{k=\infty}\widehat{\beta}_{\sigma}(x+2\pi
k)\widehat{h}(x+2\pi
k), \quad \forall x \in [0,2\pi].
\label{eqn3.3}
\end{equation}
By setting $x=2\pi u$, it follows that for $u \in [0,1],$ equation \eqref{eqn3.3} can be  modified as
\begin{equation}
\sum\limits_{k=-\infty}^{\infty}\widehat{\beta}_{\sigma}(x+2\pi k)\widehat{h}(x+2\pi k)=\sum\limits_{k=-\infty}^{\infty}\widehat{\beta}_{\sigma}(2\pi u+2\pi k)\widehat{h}(2\pi u+2\pi k).
\label{eqn3.4}
\end{equation}
Let us set
\begin{equation}
P(2\pi u):=\sum\limits_{k=-\infty}^{\infty}\widehat{\beta}_{\sigma}(2\pi u+2\pi k)\widehat{h}(2\pi u+2\pi k) , \quad\forall u \in [0,1].
\label{eqn3.5}
\end{equation}

Now \\
\begin{align}
|P(2\pi u)| & = \left|\sum\limits_{k=-\infty}^{\infty}\widehat{\beta}_{\sigma}(2\pi u+2\pi k)\widehat{h}(2\pi u+2\pi k)\right|\nonumber\\
&\geq \left|\widehat{\beta}_{\sigma}(2\pi u)\widehat{h}(2\pi u)+\widehat{\beta}_{\sigma}(2\pi u-2\pi)\widehat{h}(2\pi u-2\pi )\right |- \sum\limits_{k\in {\mathbb{Z}}\setminus \{-1,0\}}\left|\widehat{\beta}_{\sigma}(2\pi u+2\pi k)\widehat{h}(2\pi u+2\pi k)\right|.
\label{eqn3.6}
\end{align}
\noindent
We can write $h(t)$ as a  sum of the form $h(t):=h_{0}(t)+h_{1}(t)$, where $h_{0}(t):=\frac{1}{2}[h(t)+h(-t)]$ is an even function and $h_{1}(t):=\frac{1}{2}[h(t)-h(-t)]$
is an odd function. Hence
\begin{equation}
\left|\widehat{\beta}_{\sigma}(2\pi u)\widehat{h}(2\pi u)+\widehat{\beta}_{\sigma}(2\pi u-2\pi)\widehat{h}(2\pi u-2\pi )\right| > \widehat{\beta}_{\sigma}(2\pi u)\widehat{h}_{0}(2\pi u)+\widehat{\beta}_{\sigma}(2\pi (1-u))\widehat{h}_{0}(2\pi (1-u) ).
\end{equation}
 As $\widehat{\beta}_{\sigma}$ and $\widehat{h}_{0}$ are even functions,  it is sufficient to consider $u \in [0,\frac{1}{2}]$ in the above sum. Now  for  $u \in [0,\frac{1}{2}],$
\begin{align*}
\widehat{\beta}_{\sigma}(2\pi u)\widehat{h_{0}}(2\pi u) & +\widehat{\beta}_{\sigma}(2\pi (1-u))\widehat{h_{0}}(2\pi (1-u) )\\
&= \widehat{\beta}_{\sigma}(2\pi u)\int_{-l}^{l}h_{0}(t)e^{-2\pi itu}dt+\widehat{\beta}_{\sigma}(2\pi (1-u))\int_{-l}^{l}h_{0}(t)e^{-2\pi it(1-u)}dt\\
&=2\widehat{\beta}_{\sigma}(2\pi u)\int_{0}^{l}h_{0}(t)\cos(2\pi tu)dt+2\widehat{\beta}_{\sigma}(2\pi (1-u))\int_{0}^{l}h_{0}(t)\cos(2\pi t(1-u))dt\\
&\geq ||h_{0}||_{1}\left[\widehat{\beta}_{\sigma}(2\pi u)\cos(2\pi lu)+\widehat{\beta}_{\sigma}(2\pi (1-u))\cos(2\pi l(1-u))\right]\\
&\geq ||h_{0}||_{1}\left[\widehat{\beta}_{\sigma}(2\pi u)\cos(2\pi lu)+\widehat{\beta}_{\sigma}(2\pi u)\left(\frac{u}{1-u}\right)^{\sigma+1}\cos(2\pi lu)\cos(2\pi l)\right]\\
&\geq ||h_{0}||_{1}\widehat{\beta}_{\sigma}(2\pi u)\cos(2\pi lu)\left[1+\left(\frac{u}{1-u}\right)^{\sigma+1}\cos(2\pi l)\right]\\
&\geq ||h||_{1}\left(\frac{2}{\pi}\right)^{\sigma+1}\cos(\pi l)\left[1+\cos(2\pi l)\right]\\
&\geq 2||h||_{1}\left(\frac{2}{\pi}\right)^{\sigma+1}\cos^{3}(\pi l)\\
&\geq 2||h||_{1}\left(\frac{2}{\pi}\right)^{\sigma+1}\frac{1}{2^{\sigma+1}}.
\end{align*}\\
Also for $u \in [0,1],$
\begin{align*}
\sum\limits_{k\in {\mathbb{Z}}\setminus[-1,0]}|\widehat{\beta}_{\sigma}(2\pi u &+2\pi k)\widehat{h}(2\pi u+2\pi k)|\\
 & \leq ||\widehat{h}||_{\infty}\left(\sum\limits_{k=1}^{\infty}|\widehat{\beta}_{\sigma}(2\pi u+2\pi k)|+\sum\limits_{k=2}^{\infty}|\widehat{\beta}_{\sigma}(2\pi u-2\pi k)|\right)\\
&\leq ||\widehat{h}||_{\infty}\ \widehat{\beta}_{\sigma}(2\pi u)\left(\sum\limits_{k=1}^{\infty}\left|(-1)^{k(\sigma+1)}\left(\frac{u}{u+k}\right)^{\sigma+1}\right|+\sum\limits_{k=1}^{\infty}\left|(-1)^{k(\sigma+1)}\left(\frac{u}{k+1-u}\right)^{\sigma+1}\right|\right)\\
&\leq ||\widehat{h}||_{\infty}\ \widehat{\beta}_{\sigma}(2\pi u)\left(\sum\limits_{k=1}^{\infty}\left|(-1)^{k(\sigma+1)}\left(\frac{u}{u+k}\right)^{\sigma+1}\left(1+\left(\frac{u+k}{k+1-u}\right)^{\sigma+1}\right)\right|\right)\\
&\leq ||h||_{1}\left(\frac{1}{\pi}\right)^{\sigma+1}\sum\limits_{k=1}^{\infty}\left[\left(\frac{1}{k+1}\right)^{\sigma+1}+\left(\frac{1}{k}\right)^{\sigma+1}\right]\\
&\leq 2||h||_{1}\left(\frac{1}{\pi}\right)^{\sigma+1}\left(\sum\limits_{k=1}^{\infty}\left(\frac{1}{k+1}\right)^{\sigma+1}+\frac{1}{2}\right).\\
\end{align*}
\noindent
Substituting these values in \eqref{eqn3.6}, we obtain
\begin{align*}
|P(2\pi u)| &\geq \left(\left|\widehat{\beta}_{\sigma}(2\pi u)\widehat{h}(2\pi u)+\widehat{\beta}_{\sigma}(2\pi u-2\pi)\widehat{h}(2\pi u-2\pi )\right|-\sum\limits_{k\in {\mathbb{Z}}\setminus[-1,0]} \left|\widehat{\beta}_{\sigma}(2\pi u+2\pi k)\widehat{h}(2\pi u+2\pi k)\right|\right)\\
&\geq \left[2||h||_{1}\left(\frac{1}{\pi}\right)^{\sigma+1}-2||h||_{1}\left(\frac{1}{\pi}\right)^{\sigma+1}\left(\sum\limits_{k=1}^{\infty}\left(\frac{1}{k+1}\right)^{\sigma+1}+\frac{1}{2}\right)\right]\\
& = 2||h||_{1}\left(\frac{1}{\pi}\right)^{\sigma+1}\left[\frac12-\sum\limits_{k=1}^{\infty}\left(\frac{1}{k+1}\right)^{\sigma+1}\right]\\
& \geq 2||h||_{1}\left(\frac{1}{\pi}\right)^{\sigma+1} \left(\frac12 - \int_1^\infty \frac{dt}{(t+1)^{\sigma+1}}\right)\\
& = 2||h||_{1}\left(\frac{1}{\pi}\right)^{\sigma+1} \left(\frac12 - \frac{2^{-\sigma}}{\sigma}\right)> 0,
\end{align*}
for all $\sigma > 1$, since the function $\sigma\mapsto \frac12 - \frac{2^{-\sigma}}{\sigma}$ is monotonically increasing for $\sigma\geq 1$ and has value 0 at $\sigma = 1$.
\noindent
Hence, we obtain $G_{h}(z)$ has no root on the unit circle $|z|=1.$
\hfill{$\square$}
\end{pf}

\vskip 1em

The fundamental cardinal spline $L_{h,\sigma}$ of fractional order $\sigma$ with $\sigma > 1$ is an element of $L^{1}({\mathbb{R}}) \cap L^{2}({\mathbb{R}})$ because $\beta_{\sigma} \in L^{1}({\mathbb{R}}) \cap L^{2}({\mathbb{R}}),$ Furthermore, $L_{h,\sigma}$ is uniformly continuous on ${\mathbb{R}}.$

\section{Kramer's sampling theorem for local averages}

\begin{thm}
Let $\emptyset \neq I,$ $\Omega \subseteq {\mathbb{R}}$ and let $\{\phi_{k}: \mbox{  } k \in {\mathbb{Z}}\}$ be an orthonormal basis of $L^{2}(I),$ where $I$ is an interval in ${\mathbb{R}}.$ Suppose that $\{S_{k}: \mbox{  } k \in {\mathbb{Z}}\}$ is a sequence of functions $S_{k}:\Omega\rightarrow {\mathbb{C}}$ and ${\bf t}:=\{t_{k} \in {\mathbb{R}}: \mbox{  } k \in {\mathbb{Z}}\}$ a numerical sequence in $\Omega$ and the averaging function $h(t)$ is compactly supported in $L^{1}({\mathbb{R}})$ satisfying the conditions
\begin{itemize}
\item[\emph{D1.}] $S_{k}\star h(t_{l})=a_{k}\delta_{kl},$ $(k,l) \in {\mathbb{Z}} \times {\mathbb{Z}},$ where $a_{k} \neq 0$;
\item[\emph{D2.}] $\displaystyle{\sum\limits_{k \in {\mathbb{Z}}}|S_{k}(t)|^{2}< \infty},$ for each $t \in \Omega.$
\end{itemize}
Define a function $K:I\times \Omega \rightarrow \mathbb{C}$ by
\begin{equation}
K(x,t):=\sum\limits_{k \in {\mathbb{Z}}}S_{k}(t)\overline{\phi_{k}}(x),
\label{eqnk1}
\end{equation}
and a linear integral transform $T$ on $L^{2}(I)$ by
\begin{equation}
(TF)(t):= \int_{I}F(x)K(x,t)dx.
\label{eqnk2}
\end{equation}
Then $T$ is well defined and injective. Furthermore, if the range of $T$ is denoted by
$$\mathcal{H}:=\{f: \mbox{  } {\mathbb{R}}\rightarrow {\mathbb{C}}: \mbox{  } f(t)= \int_{I}F(x)K(x,t)dx, \mbox{  } F \in L^{2}(I)\},$$
then
\begin{enumerate}
\item[\emph{(i)}]
$(\mathcal{H}, \langle \cdot, \cdot \rangle _{\mathcal{H}})$ is a Hilbert space isometrically isomorphic to $L^{2}(I),$ i.e, $\mathcal{H}\cong L^{2}(I),$ when endowed with the inner product
$$\langle f,g\rangle_{\mathcal{H}}:= \langle F,G \rangle_{L^{2}(I)},$$
where $f(t)=TF(t)= \int_{I}F(x)K(x,t)dx$ and $g(t)=TG(t)= \int_{I}G(x)K(x,t)dx.$
\item[\emph{(ii)}]
$\{S_{k}: \mbox{  } k \in {\mathbb{Z}}\}$ is an orthonormal basis for $\mathcal{H}.$
\item[\emph{(iii)}]
Each function $ f \in \mathcal{H}$ can be recovered from its samples on the sequence $\{t_{k}: k \in {\mathbb{Z}}\}$ via the formula
$$\displaystyle{f(t)= \sum\limits_{k \in {\mathbb{Z}}}f\star h(t_{k}) \frac{S_{k}(t)}{a_{k}}}.$$
The above series converges absolutely and uniformly on subsets of ${\mathbb{R}},$ where $||K(\cdot,t)||_{L^{2}(I)}$ is bounded.
\end{enumerate}
\label{averagek}
\end{thm}

\begin{pf}
By the Cauchy-Schwartz inequality, the linear integral transform \ref{eqnk2} is well defined for each $t \in \Omega,$ since $F$ and $K(\cdot, t)$ are in $L^{2}(I).$ Now,
\begin{eqnarray*}
K \star h (x,t_{k})&=&\sum\limits_{n \in {\mathbb{Z}}}(S_{n}\star h)(t_{k})\overline{\phi_{n}}(x)\\
&=&a_{k}\overline{\phi_{k}}(x)
\end{eqnarray*}
Further, the transformation \ref{eqnk2} is one to one because $\left\{K \star h (x,t_{k})=a_{k}\overline{\phi_{k}}(x)\right\}_{k=1}^{\infty}$ is a complete orthogonal sequence for $L^{2}(I).$

Let $\mathcal{H}$ be the range of the integral transform  endowed with the norm $||f||_{\mathcal{H}}=||F||_{L^{2}(I)},$ where $f=T(F).$
Consider $f(t)= \int_{I}F(x)K(x,t)dx$ and $g(t)= \int_{I}G(x)K(x,t)dx.$\\
Using the polarization identity, we obtain
\begin{eqnarray*}
\langle f,g\rangle_{\mathcal{H}}&=&\frac{1}{4}\left[4||f||_{\mathcal{H}} ||g||_{\mathcal{H}}\right]\\
&=&||F||_{L^{2}(I)}||G||_{L^{2}(I)}\\
&=&\langle F,G \rangle_{L^{2}(I)}.
\end{eqnarray*}
Therefore, $(\mathcal{H}, \langle \cdot, \cdot \rangle _{\mathcal{H}})$ is a Hilbert space isometrically isomorphic to $L^{2}(I),$ with the inner product
\begin{equation}
\langle f,g\rangle_{\mathcal{H}}= \langle F,G \rangle_{L^{2}(I)},
\label{eqnk3}
\end{equation}
Now we prove $\{S_{n}: \mbox{  } n \in {\mathbb{Z}}\}$ is an orthonormal basis for $\mathcal{H}.$\\
For every $k$,
$T(\phi_{k})=S_{k}$ and hence  we obtain
$$ \langle S_{n}, S_{k} \rangle_{\mathcal{H}}= \langle \phi_{n}, \phi_{k} \rangle_{L^{2}(I)}.$$
Therefore $\{S_{n}: \mbox{  } n \in {\mathbb{Z}}\}$ is an orthonormal basis for $\mathcal{H}.$\\
\vskip 1em

\indent Expanding the functions $f \in {\mathcal{H}}$ with respect to the orthonormal basis $\{S_{n}(t)\}_{n=-\infty}^{\infty},$ we have
\begin{equation}
f(t)= \sum\limits_{n=-\infty}^{\infty} \langle f, S_{n} \rangle_{\mathcal{H}}S_{n}(t),
\label{eqnk4}
\end{equation}
where the convergence is in the ${\mathcal{H}}$ norm sense and hence pointwise in $\Omega.$
By (i), the isometry between ${\mathcal{H}}$ and $L^{2}(I),$ we obtain
\begin{eqnarray}
\nonumber
\langle f, S_{n} \rangle_{\mathcal{H}}&=& \langle F, \phi_{n} \rangle_{L^{2}(I)}\\
\label{eqnk5}
&=& \int_{I} F(x)\overline{\phi_{n}}(x)dx,
\end{eqnarray}
where $T(F)=f.$
Using the integral transform \ref{eqnk2},
\begin{eqnarray*}
\nonumber
f \star h(t_{n})&=& \int_{I} F(x)K \star h(x, t_{n})dx\\
\nonumber
&=&\int_{I} F(x)a_{n}\overline{\phi_{n}}(x) dx.\\
\end{eqnarray*}
Hence
\begin{equation}
\label{eqnk6}
\frac{f \star h(t_{n})}{a_{n}}=\int_{I} F(x) \overline{\phi_{n}}(x) dx.
\end{equation}
By \ref{eqnk5} and \ref{eqnk6}, we obtain
$$\langle f, S_{n} \rangle_{\mathcal{H}}=\frac{f \star h(t_{n})}{a_{n}}.$$
Therefore by \ref{eqnk4}
\begin{eqnarray*}
f(t)&=& \sum\limits_{n=-\infty}^{\infty} \langle f, S_{n} \rangle_{\mathcal{H}}S_{n}(t)\\
&=&\sum\limits_{n=-\infty}^{\infty} f \star h(t_{n})\frac{S_{n}(t)}{a_{n}}.
\end{eqnarray*}
The above series converges absolutely and uniformly on subsets of ${\mathbb{R}},$ where $||K(\cdot,t)||_{L^{2}(I)}$ is bounded.
\hfill{$\square$}
\end{pf}

In theorem \ref{averagek}, we choose $\Omega:={\mathbb{R}},$ ${\bf t}:= {\mathbb{Z}},$ $a_{k}=1$, for all $k \in {\mathbb{Z}},$ and for the interpolating function $S_{k}=L_{h,\sigma}(\cdot -k),$ $k \in {\mathbb{Z}}.$ Then we obtain the average sampling theorem for fundamental splines of fractional order.
\begin{thm}
Let $\emptyset \neq I \subseteq {\mathbb{R}}$ and let $\{\phi_{k}: \mbox{  } k \in {\mathbb{Z}}\}$ be an orthonormal basis of $L^{2}(I).$ Let $L_{h,\sigma}$ denote the fundamental cardinal spline of fractional order $\sigma \in T_{h}.$ Then the following hold:
\begin{enumerate}
\item[\emph{(i)}]
The family $\{L_{h,\sigma}(\cdot -k): k \in {\mathbb{Z}}\}$ is an orthonormal basis of the Hilbert space $(\mathcal{H},\langle \cdot,\cdot \rangle_{\mathcal{H}}),$ where $\mathcal{H} = T(L^{2}(I))$ and $T$ is  the injective integral operator
$$\displaystyle{T(F)=\sum\limits_{k \in {\mathbb{Z}}}\langle F,\phi_{k}\rangle_{L^{2}(I)}L_{h,\sigma}(\cdot -k)}, \mbox{     } F\in L^{2}(I).$$
\item[\emph{(ii)}]
Every function $f \in \mathcal{H}\cong L^{2}(I)$ can be recovered from its samples on the integers via
\begin{equation}
f(\cdot)=\sum\limits_{k \in {\mathbb{Z}}}f\star h(k)L_{h,\sigma}(\cdot -k),
\label{kramaeq}
\end{equation}
where the series (\ref{kramaeq}) converges absolutely and uniformly on all subsets of ${\mathbb{R}}.$
\end{enumerate}
\end{thm}
\begin{pf}
Condition D1. for $S_{k}=L_{h,\sigma}(\cdot -k),$ $k \in {\mathbb{Z}}$, in Theorem \ref{averagek} can be modified as
$$L_{h,\sigma}\star h(n)=\delta_{n,0}, \mbox{   } n \in {\mathbb{Z}}.$$
This condition is established by the equation (\ref{finterh}). The condition D2. is also verified because the fundamental cardinal spline $L_{h,\sigma}$ is an element of $L^{1}({\mathbb{R}}) \cap L^{2}({\mathbb{R}}).$ Since the unfiltered splines $\{\beta_{\sigma}(\cdot -k): k \in {\mathbb{Z}}\}$ already form a Riesz basis of the space $V$ (see, \cite{uns2}), $||K(\cdot,t)||_{L^{2}(I)}$ is bounded on ${\mathbb{R}}.$
Hence by Theorem \ref{averagek} the statements (i) and (ii) hold.
\hfill{$\square$}

\end{pf}

\begin{rmk}
In \cite{han},  D. Han et al. investigated reproducing kernel Hilbert spaces on a set $\Omega$ which contains a given countable subset $\Lambda\subset \Omega$ as a sampling set. A similar analysis may be performed for Kramer-type samplings as well. Likewise, the results obtained in \cite{sun3,sun4,vet} on sampling with finite rates of innovation can also be applied to fractional spline spaces. These questions will be investigated in a forthcoming paper.
\end{rmk}

\noindent {\bf Acknowledgement:-} The third author would like to thank the management of Sri Sivasubramaniya Nadar College of Engineering. The second author was partially supported by DFG grant MA 5801/2-1.



\begin{thebibliography}{99}
\bibitem{ald} A. Aldroubi and M. Unser. Sampling procedure in function spaces and asymptotic equivalence with Shannon's sampling theory, Numerical Functional Analysis and Optimization, 15(1) (1994), 1--21.
\bibitem{ald1} A. Aldroubi and K.-H. Gr\"ochenig. Beurling-Landau-type theorems for nonuniform sampling in shift-invariant spaces, J. Fourier Anal. Appl. 6(1) (2000), 91--101.
\bibitem{ald2} A. Aldroubi and K.-H. Gr\"ochenig. Nonuniform sampling and reconstruction in shift-invariant spaces, SIAM Review 43(4) (2001), 585--562.
\bibitem{ald3} A. Aldroubi and M. Unser. Sampling procedure in function spaces and asymptotic equivalence with Shannon's sampling theory, Numerical Functional Analysis and Optimization, 15(1) (1994), 1--21.
\bibitem{uns} M. Unser. Splines: A perfect fit for signal and image processing, IEEE Signal Processing Magazine 16 (1999), 22--38.
\bibitem{uns1} M. Unser. Sampling 50 years after Shannon, Proceeding IEEE 88 (2000), 569--587.
\bibitem{gro} K.-H. Gr\"{o}chenig. Reconstruction algorithms in irregular samplings, Mathematics of Computation 59 (1992), 181--194.
\bibitem{buz} P.L. Butzer, W. Engels, S. Ries and R.L. Stens. The Shannon sampling series and the reconstruction of signals in terms of linear, quadratic and cubic splines, SIAM J. Appl. Math., 46(2) (1986), 299--323.
\bibitem{sun} W. Sun and X. Zhou. Average sampling in spline subspaces, Applied Mathematics Letters, 15 (2002), 233--237.
\bibitem{sun1} W. Sun and X. Zhou. Reconstruction of functions in spline subspaces from local averages, Proc. Amer. Math. Soc., 131 (2003), 2561--2571.
\bibitem{chu} C. K. Chui. An Introduction to Wavelets, Academic Press,  1992.
\bibitem{sch} I.J. Schoenberg. Cardinal Spline Interpolation, SIAM Regional Conference Series in Applied Mathematics, 1973.
\bibitem{uns2} M. Unser and Th. Blu. Fractional splines and wavelets, SIAM Review., 42(1) (2000), 43--67.
\bibitem{for} B. Forster and P. Massopust. Interpolation with fundamental splines of fractional order, Proceedings of SampTA 2011, 1--4.
\bibitem{deb} C. de Boor, K. H\"ollig and S. Riemenschneider. Bivariate cardinal interpolation by splines on a three-direction mesh, Ill. J. Math., 29(4) (1985) 533--566.
\bibitem{chu1} C.K. Chui, K. Jetter and J.D. Ward. Cardinal interpolation by multivariate splines, Mathematics of Computation, 48(178) (1987), 711--724.
\bibitem{gar1} A. G. Garcia. Orthogonal sampling formulas: A unified approach, SIAM Review, 42(3) (2000), 499--512.
\bibitem{kra} H. P. Kramer. A generalized sampling theorem, J. Math. Phys., 63 (1957), 68--72.
\bibitem{zho} X. Zhou and W. Sun. On the sampling theorem for wavelet subspaces, J. Fourier Anal. Appl., 5(4) (1999), 347--354.
\bibitem{sun2} W. Sun and X. Zhou. Frames and sampling theorem, Science in China (Series A), 41(6) (1998), 606--612.
\bibitem{baer} H. Batemann and A. Erd{\'e}lyi. Higher Transcendental Functions I, McGraw Hill Book Company, 1953.
\bibitem{han} D. Han, M. Z. Nashed and Q. Sun. Sampling expansions in reproducing kernel Hilbert and Banach spaces, Numerical Functional Analysis and Optimization, 30(9-10) (2009), 971--987.
\bibitem{Ald4} A. Aldroubi, Q. Sun and W.-S. Tang. Nonuniform average sampling and reconstruction in multiply generated shift-invariant spaces, Constr. Approx., 20 (2004), 173--189.
\bibitem{Ald5} A. Aldroubi, Q. Sun and W.-S. Tang. Convolution, average sampling, and a Calder\'on resolution of the identity for shift-invariant spaces, J. Fourier Anal. Appl., 11 (2005), 215--244.
\bibitem{sun3} Q. Sun, Nonuniform Average Sampling and Reconstruction of signals with finite rate of innovation, SIAM J. Math. Anal., 38 (2006),1389--1422.
\bibitem{sun4} Q. Sun. Frames in spaces with finite rate of innovation, Adv. Comput. Math., 28 (2008), 301--329.
\bibitem{vet} M. Vetterli, P. Marziliano and Th. Blu. Sampling signals with finite rate of innovation, IEEE Trans. Signal Process., 50 (2002), 1417--1428.
\end{thebibliography}
\end{document}